\documentclass{amsart}
\usepackage{amssymb,amsfonts,amsmath}
\usepackage{graphicx}

 \newtheorem{thm}{Theorem}[section]
 
 \newtheorem{lem}[thm]{Lemma}
 \newtheorem{prop}[thm]{Proposition}
 \theoremstyle{definition}
 
 \theoremstyle{remark}
 \newtheorem{rem}[thm]{Remark}
 
 \numberwithin{equation}{section}

\newcommand{\eps}{\varepsilon}

\newcommand{\rank}{{\rm rank\, }}
\newcommand{\Toep}{{\rm Toep\, }}

\newcommand{\seq}[1]{\left\langle#1\right\rangle}
\newcommand{\set}[1]{\left\{#1\right\}}

\begin{document}
\title
{Eigenvalues of rank one perturbations of unstructured matrices}
\author[A.C.M. Ran]{Andr\'e C.M. Ran}
\address{Afdeling Wiskunde, Faculteit der Exacte Wetenschappen, Vrije
Universiteit Amsterdam, De Boelelaan 1081a, 1081 HV Amsterdam, The Netherlands}

\author[M. Wojtylak]{ Micha\l{} Wojtylak}
\address{Instytut Matematyki, Wydzia\l{} Matematyki i Informatyki, Uniwersytet Jagiello\'nski, Krak\'ow, ul. \L ojasiewicza 6, 30-348 Krak\'ow, Poland}

\thanks{The second author gratefully acknowledges the assistance of the Polish Minstery of Higher Education and Science grant NN201 546438. }
\maketitle

\begin{abstract}
Let $A$ be a fixed complex matrix and let $u,v$ be two vectors. The eigenvalues of matrices
$A+\tau uv^\top $  $(\tau\in\mathbb{R})$ form a system of intersecting curves. The dependence
 of the intersections on the vectors $u,v$ is studied.
\end{abstract}

\section*{Introduction}

\begin{figure}[htb]
\includegraphics[width=12cm]{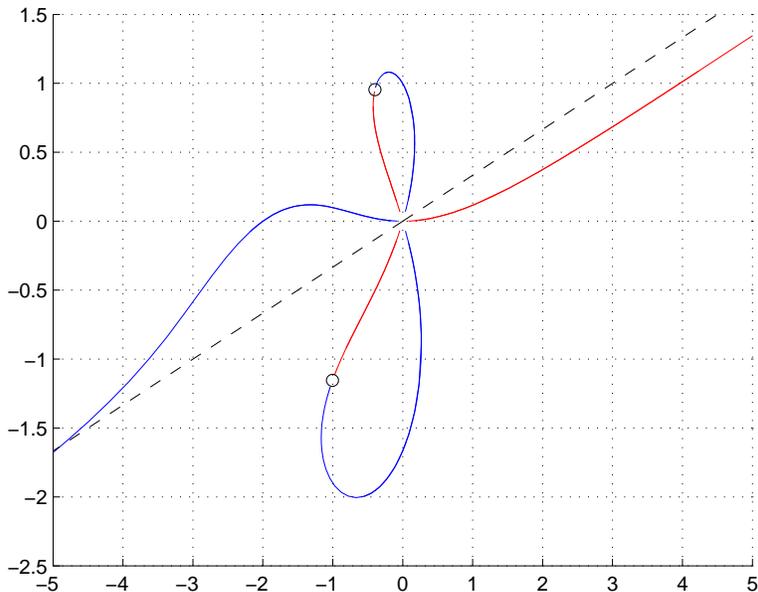}
 \caption{Eigenvalues of $B(\tau)=J_3(0)+\tau uv^\top $ for $\tau>0$ in red, for $\tau<0$ in blue, see also Remark \ref{desc}.}\label{nodoublevlees}
\end{figure}

The motivation for this paper is the following numerical experiment. Take a  matrix
$A\in\mathbb{C}^{n\times n}$ and nonzero vectors $u,v\in\mathbb{C}^n$ and plot the set
\begin{equation}\label{problemstatement}
\{\sigma(A+\tau uv^\top ):\tau\in\mathbb{R}\}.
\end{equation}
It is well known that  above set consists of a finite number of curves, that intersect only in
a finite number of points. However, it appears that for $u,v\in\mathbb{C}^n$ chosen randomly
from a continuous distribution on $\mathbb{C}^n$ there are no intersection points except,
possibly, the of spectrum of $A$. Furthermore, for all $\tau\in\mathbb{R}\setminus\{0\}$ all
eigenvalues of $A+\tau uv^\top $, that are not eigenvalue of $A$, are simple. A typical case
for $A=J_3(0)$ is shown of Figure \ref{nodoublevlees}, note that the only intersection of the
eigenvalue curves is at $0\in\sigma(A)$. Since it appears that the intersection points outside
$\sigma(A)$ are multiple eigenvalues of $A+\tau uv^\top$ (cf. Proposition \ref{leigs}(ii)), we
will be also interested in a problem of existence of multiple eigenvalues of $A+\tau uv^\top$
for some $\tau\in\mathbb{C}$.

Some light on the phenomenon of lack of double eigenvalues in the numerical simulations is put
by the following marvelous result of H\"ormander and Melin \cite{HM}. Let the Jordan canonical
form of the matrix $A$ be
$$
A\cong\bigoplus_{j=1}^r \bigoplus_{i=1}^{k_j} J_{n_{j,i}} (\lambda_j) ,
$$
where the Jordan blocks $J_{n_{j,i}} (\lambda_j)$ corresponding to each eigenvalue $\lambda_j$
($j=1,\dots r$) are in decreasing order, i.e. $n_{j,1}\geq n_{j,2}\geq\dots\geq n_{j,k_j}$.
Then for generic $u$ and $v$ (i.e. for all $u$ and $v$ except a `small' set, see
Preliminaries) the Jordan form of $A+uv^\top $ is the following
$$
A+uv^\top \cong\bigoplus_{j=1}^r \bigoplus_{i=2}^{k_j} J_{n_{j,i}} (\lambda_j)\oplus
\bigoplus_{h=1}^l J_1{\mu_h},
$$
where $\mu_h\neq \mu_{h'}$ for $h\neq h'$. In other words, for each eigenvalue $\lambda_j$
($j=1,\dots r$) only the largest chain in the Jordan structure is destroyed and there appears
a structure of simple eigenvalues instead.

The behavior of  eigenvalues of $A+\tau uv^\top $ as functions of $\tau$ for small values of
$\tau$ is also well known, see, e.g., \cite{Lidskii, Vishik-Lyusternik} and
\cite{Baumgartel,Kato,Moro-Burke-Overton,Sa1,Wilkinson}. Namely, for small values of $|\tau |$
and for generic $u$ and $v$ for each $j=1,\dots r$ there are $n_{j,1}$ simple eigenvalues
$\mu_{j,k}(\tau)$, $k=1, \cdots , n_{j,1}$ of $A+\tau uv^\top $ in a punctured neighborhood of
$\lambda_j$, and they are given by
\begin{equation}\label{asymptoticssmallt2}
\mu_{j,k}(\tau) =\lambda_j + \tau^{1/{n_{j,1}}}\cdot \left( c_j \right)^{1/n_{j,1}} \cdot
\exp\left(\frac{-2\pi i k }{n_{j,1}}\right) + O(\tau^{2/n_{j,1}}) ,
\end{equation}
where the number $c_j$ can be expressed explicitly in terms of $A,u$ and $v$; see \cite{Sa1},
Proposition 1. That is, the eigenvalues $\mu_{j,k}(\tau)$ are approximately given by the roots
of the polynomial equation
\begin{equation}\label{asymptoticssmallt}
(\mu-\lambda_j)^{n_{j,1}}= \tau \cdot c_j ,\quad j=1,\dots r.
\end{equation}
However, neither the H\"ormander--Mellin result nor the above small $\tau$ asymptotic of
eigenvalues does not explain the lack of crossing of eigenvalue curves that appears in
numerical simulations.  The purpose of the present paper is to show that this behavior is
indeed `generic' although the notion of genericity will have some different shades.

For historical reasons let us mention two works prior to the H\"ormander--Mellin paper, in
\cite{Thompson} the invariant factors of a one--dimensional perturbation are considered and in
\cite{Krupnik} the perturbation theory for normal matrices is developed. The result by
H\"ormander--Mellin lay dormant for about a decade before being rediscovered independently by
Dopico and Moro \cite{DM} and Savchenko \cite{Sa1,Sa2}. Since that time the interest in topic
has grown up, see e.g. \cite{MMRR1,MMRR2} for an alternative proof using ideas from systems
theory and for perturbation theory for structured matrices. Although the results presented
below concern a similar matter the reasonings are independent of the previous work and the
content of the paper is self--contained. The main outcome  are Theorems \ref{thm:generict},
\ref{thm:largeasymp}, \ref{thm:2}, \ref{prop:3}  and \ref{propC}. First four of them allow the
parameter $\tau$ to be complex, while in the last one we return to the real parameter $\tau$.
This collection gives a complete description of the generic behavior of the set in
\eqref{problemstatement}.

\section{Preliminaries}

In this section, we gather some known results which will be the basis for our further
investigation. An important technique used in this paper is the resultant. Let
$$
q_1(\lambda)=a_{n_1}\lambda^{n_1}+\cdots +a_0,\qquad q_2(\lambda)=b_{n_2}\lambda^{n_2}+\cdots
+b_0
$$
be two complex polynomials. By $S(q_1,q_2)$ we denote the \textit{Sylwester resultant matrix}
of $q_1$ and $q_2$:
 \begin{equation}\label{Sylwester} S(q_1,q_2)=\begin{bmatrix}
a_{n_1}&  \cdots & a_0 & 0 &  \cdots &\cdots & 0 \\
0 & a_{n_1} &   \cdots &  a_0  & 0 & \cdots & 0 \\
\vdots & & \ddots &  & \ddots&  & 0\\
0 & \cdots &  \cdots& 0 &  a_{n_1} &  \cdots &  a_0\\
b_{n_2}&  \cdots & b_0 & 0 &  \cdots &\cdots & 0 \\
0 & b_{n_2}&   \cdots &  b_0  & 0 & \cdots & 0 \\
\vdots & & \ddots &  & \ddots&  & 0\\
0 & \cdots &  \cdots& 0 &  b_{n_2} &  \cdots &  b_0\\
\end{bmatrix}\in\mathbb{C}^{(n_1+n_2)\times(n_1+n_2)}.
\end{equation}
It is well known that $q_1$ and $q_2$ have a common root if and only if $\det S(q_1,q_2)=0$.

Let $A\in\mathbb{C}^{n\times n}$ and let $u,v\in\mathbb{C}^n$. Occasionally we will use the
notation
$$
B(\tau )= A+\tau uv^\top,\quad \tau\in\mathbb{C},
$$
remembering, nevertheless, that we are interested in the $(u,v)$--dependence of the spectral
structure of $B(\tau)$. Recall that an eigenvalue $\lambda_0$ of $B\in\mathbb{C}^{n\times n}$
is called  \emph{non--derogatory\/}  if $\dim \ker (B-\lambda_0)=1$. The following result may
be found in \cite{Sa1}, Lemma 5, for completeness sake we include a proof.

\begin{lem}\label{prop:1}
Let $A\in\mathbb{C}^{n\times n}$  and let $u,v\in\mathbb{C}^n$. Then for all
$\tau\in\mathbb{C}\setminus \set0$ all eigenvalues of $B(\tau)$ that are not eigenvalues of
$A$
 are non-derogatory.
\end{lem}

\begin{proof}
Let $\lambda_0\in\sigma(B(\tau))\setminus \sigma(A)$ and let $\tau\neq0$. Using  the fact that
$\rank(X+Y) \leq \rank X + \rank Y$ for any compatible matrices $X,Y$ we obtain
\begin{align*}
n=\rank (A-\lambda_0 ) &\leq \rank (A+\tau uv^\top   -\lambda_0 )+\rank (\tau uv^\top  ) =\\
&= \rank (B(\tau)-\lambda_0 ) +1 ,
\end{align*}
which shows that $\rank(B(\tau)-\lambda_0  ) \geq n-1$. Hence $\dim \ker (B(\tau)-\lambda_0) =
1$ and so $\lambda_0$ is a non-derogatory eigenvalue of $B(\tau)$.
\end{proof}

Following \cite{MMRR1} we say that a subset $\Omega$ of $\mathbb C^n$ is \textit{generic} if
$\Omega$ is not empty and the complement $\mathbb{C}^{n} \setminus\Omega$ is contained in a
(complex) algebraic set which is not $\mathbb{C}^n$. In such case $\mathbb C^n\setminus\Omega$
is nowhere dense and of $2n$--dimensional Lebesgue measure zero. We use the phrase \textit{for
generic $v\in \mathbb C^n$} as an abbreviation of: `there exist
a generic $\Omega\subseteq\mathbb C^n$ such that for all $v\in\Omega$'. 
Our main results, except Theorem \ref{propC}, have the following form:
 \begin{itemize}
 \item[] \textit{Let $A\in\mathbb{C}^{n\times n}$. Then for generic $u$ and $v$ ...},
 \end{itemize}
which should be read formally as
  \begin{itemize}
 \item[] \textit{For every $A\in\mathbb{C}^{n\times n}$ there exists a
generic subset $\Omega$ of $\mathbb{C}^{2n}$, possibly dependent on $A$,  such that for
$(u,v)\in \Omega$...}.
 \end{itemize}

Most of our reasoning are independent of a choice of basis. Let $T$ be an invertible matrix. Then
$$
 T(A+\tau uv^\top )T^{-1}=TAT^{-1} + \tau (Tu)(v^\top T^{-1}).
$$
In consequence, the Jordan structures of the matrices $A+\tau uv^\top $ and $TAT^{-1} + \tau
(Tu)(v^\top T^{-1})$ are identical. In other words the transformation
\begin{equation}\label{matrixtransform}
(A,u,v^\top )\mapsto (TAT^{-1},Tu,v^\top T^{-1})
\end{equation}
preserves the spectral structure of $B(\tau)$ for all $\tau\in\mathbb{R}$. Let  $T_A$ be the
transformation of $A$ to its Jordan canonical form, that is
\begin{equation}\label{JordanA}
A'=T_AAT_A^{-1}=\bigoplus_{j=1}^r \bigoplus_{i=1}^{k_j} J_{n_{j,i}} (\lambda_j),
\end{equation}
where $J_k(\lambda)$ denotes the Jordan block of size $k$ with the diagonal entries equal
$\lambda$ and the entries on the first \underbar{upper}--diagonal equal one and
\begin{equation}\label{n's}
n_{j,1}\geq n_{j,2}\geq\cdots \geq n_{j,k_j},\quad j=1,\dots, r.
\end{equation}
We will describe now a special instance of the transformation $T$ that consists of two steps,
i.e. $T=T_vT_A$. Let $T_A$ be as above, next we decompose $u'=T_Au$ and $v'{}^\top =v^\top
T_A^\top $ according to the Jordan form of $A'$ as follows:
\begin{equation}\label{uform}
u'=\begin{bmatrix} u'_{1} \\ u'_{2} \\ \vdots \\ u'_{p_u}\end{bmatrix}, \qquad u'_j=\begin{bmatrix} u'_{j,1} \\
u'_{j,2} \\ \vdots \\ u'_{j,k_j}\end{bmatrix}, \qquad u'_{j,i}= \begin{bmatrix} u'_{j,i,1} \\
u'_{j,i,2}\\ \vdots
\\ u'_{j,i,n_{j,i}}\end{bmatrix}\in\mathbb{C}^{n_{j,i}} ,
\end{equation}
and
\begin{equation}\label{vform1}
v'=\begin{bmatrix} v'_{1} \\ v'_{2} \\ \vdots \\ v'_{p_u}\end{bmatrix}, \qquad v'_j=\begin{bmatrix} v'_{j,1} \\
v'_{j,2} \\ \vdots \\ v'_{j,k_j}\end{bmatrix}, \qquad v'_{j,i}= \begin{bmatrix} v'_{j,i,1} \\
v'_{j,i,2}\\ \vdots
\\ v'_{j,i,n_{j,i}} \end{bmatrix}\in\mathbb{C}^{n_{j,i}}.
\end{equation}
We put
$$
T_v=\bigoplus_{j=1}^r \bigoplus_{i=1}^{k_j} \Toep(v_{j,i}) ,
$$
where by $\Toep(w)$  we denote the $k\times k$ upper-triangular Toeplitz matrix whose first
row is given by $w\in\mathbb{C}^k$. Obviously $T_v$ commutes with $A$. Now note that for
generic $v$ one has
\begin{equation}\label{generic}
v'_{j,i,1}\neq 0\qquad  i=1,\dots, k_j,\ j=1,\dots,r,
\end{equation}
which implies that $T_v$ is invertible, consequently $T_vA'T_v^{-1} =A'$. Furthermore,
$v''{}^\top =v'{}^\top  T_v^{-1}$ has the following form
\begin{equation}\label{vform}
v''_{j,i}= \begin{bmatrix} 1 \\ 0 \\ \vdots \\ 0 \end{bmatrix}\qquad  i=1,\dots, k_j,\
j=1,\dots,r.
\end{equation}
The triplet $(TAT^{-1}, Tu, v^\top T^{-1})$, where $T=T_vT_A$, will be called  the
\textit{Brunovsky form of $(A,u,v^\top )$}, cf. \cite{Bru}. Note the following simple lemma,
that will allow us to reduce the problem of genericity in $u$ and $v$ to a problem of
genericity in $u$ with a fixed $v$.
\begin{lem}\label{l:Brun}
If $\Omega_0$ is a generic subset of $\mathbb{C}^n$ then the set
 $$
\{(u,v)\in\mathbb{C}^{2n}:T_v\text{ is invertible }, T_vT_Au\in\Omega_0\}
 $$
is a generic subset of $\mathbb{C}^{2n}$.
\end{lem}


\section{The characteristic polynomial of $B(\tau)$}

The present  section contains the basic tools used in the paper. Namely, we introduce the
polynomial $p_{uv}$ and provide a formula for the characteristic polynomial of $B(\tau)$.

The minimal polynomial of $A$ will be denoted by $m(\lambda)$. Everywhere in the paper
\eqref{JordanA} and \eqref{n's} are silently assumed, consequently one has
\begin{equation}\label{m}
m(\lambda)= \prod_{j=1}^r (\lambda-\lambda_j)^{n_{j,1}}.
\end{equation}
 We also put
\begin{equation}\label{formulap}
p_{uv}(\lambda)=  m(\lambda)\cdot v^\top  (\lambda -A)^{-1} u .
\end{equation}
Note that $p_{uv}$ is invariant on the transformation \eqref{matrixtransform}. Transforming
$A$ to it Jordan form we easily see that $p_{uv}$ is a polynomial of degree at most $\deg
m-1$. The following Lemma plays an essential role in the further reasoning.

\begin{lem}\label{pdegl-1}
For generic $u$ and $v$ the polynomial $p_{uv}$ is  of degree $\deg m-1$ and has no double
roots and no common roots with $m$.
\end{lem}

\begin{proof}
Using Lemma \ref{l:Brun} and the fact that $p_{uv}$ is invariant on the transformation
\eqref{matrixtransform} we may assume that $A$ is in the Brunovsky canonical form and treat
$v$ as fixed. For simplicity consider the case when $A$ consists of one Jordan block only,
i.e.
$$
A=J_n(\lambda_1),\quad v=\begin{bmatrix} 1\\0\\\vdots\\ 0\end{bmatrix},\quad u=\begin{bmatrix}
u_1\\u_2\\\vdots\\ u_n\end{bmatrix}.
$$
Then $m(\lambda)=(\lambda-\lambda_1)^n$ and
$$
(\lambda-A)^{-1}=\Toep([(\lambda-\lambda_1)^{-1},(\lambda-\lambda_1)^{-2},\dots,
(\lambda-\lambda_1)^{-n}]^\top).
$$
Consequently,
$$
p_{uv}(\lambda)=u_1(\lambda-\lambda_1)^{n-1}+\cdots +u_{n-1}(\lambda-\lambda_1)+u_n.
$$
Hence, the generic assumption $u_1\neq 0$ implies that $\deg p_{uv} =\deg m-1$. Further on,
the generic assumption $u_n\neq 0$ implies that $p_u$ and $m$ do not have common roots. To
prove that for generic $u$ the polynomial $p_{uv}$ has simple roots only  let us consider the
Sylwester resultant matrix $S(p_{uv},p_{uv}')$. Note that $\det S(p_{uv},p_{uv}')$ is a
nonzero polynomial in $u$. Hence, the equation $\det S(p_{uv},p_{uv}')=0$ defines a proper
algebraic subset of $\mathbb{C}^n$.

 The general case follows by similar arguments from the equation
$$
p_{uv}(\lambda)=m(\lambda)\cdot\sum_{j=1}^r\sum_{i=1}^{k_j} v_{j,i}^\top (\lambda-
J_{n_{j,i}}(\lambda_j))^{-1}u_{j,i}.
$$
\end{proof}

We put
$$
q(\lambda)=\prod_{i=1}^r\prod_{j=2}^{k_i} (\lambda-\lambda_i)^{n_{i,j}}
=\frac{\det(\lambda-A)}{m(\lambda)},
$$
with the convention $\prod_{2}^1:=1$. We also define the family of polynomials $p_{uv,\tau}$
by
\begin{equation}\label{ptaudef}
p_{uv,\tau}(\lambda)=m(\lambda)-\tau p_{uv}(\lambda),\quad \tau\in\mathbb{R}.
 \end{equation}

\begin{prop}\label{leigs}
Let $A\in\mathbb{C}^{n\times n}$, then the following statements hold.
 \begin{itemize}
 \item[(i)] For every $u,v\in\mathbb{C}^n$, $\tau\in\mathbb{C}$  the characteristic polynomial of $A+\tau uv^\top
$ equals $q\cdot p_{uv,\tau}$.
\item[(ii)] For every $u,v\in\mathbb{C}^n$, $\tau_1,\tau_2\in\mathbb{C}$ with
$\tau_1\neq\tau_2$ one has
$$
\sigma(A+\tau_1 uv^\top)\cap\sigma(A+\tau_2 uv^\top)\subseteq\sigma(A).
$$
 \item[(iii)] For generic $u$ and $v$ and all
$\tau\in\mathbb{C}\setminus\set0$ there are exactly $\deg m$, counting algebraic
multiplicities, eigenvalues of $A+\tau uv^\top $ that are not eigenvalues of $A$.

\end{itemize}
\end{prop}

Point (iii)  shows that the only crossings of the eigenvalue curves in
\eqref{problemstatement} are the multiple eigenvalues of $A+\tau uv^\top$ for some
$\tau\in\mathbb{R}$.

\begin{proof}
(i) For any $u,v\in\mathbb{C}^n$, $\tau\in\mathbb{C}$ and
$\lambda\in\mathbb{C}\setminus\sigma(A)$ we have (cf. \cite{Sa1}, Lemma 1)
\begin{align}\label{handycomputation}
 \det (\lambda - (A+\tau uv^\top  ))&= \det\left( (\lambda-A)(I-(\lambda
-A)^{-1}\tau uv^\top  ) \right)\nonumber
\\
&= \det(\lambda-A) \det(I-(\lambda -A)^{-1}\tau uv^\top  )\nonumber \\
&=
\det(\lambda-A) ( 1- \tau v^\top  (\lambda -A)^{-1} u)\nonumber
\end{align}
Dividing both sides by $q$ and emploing \eqref{formulap} we obtain
\begin{equation}\label{nicepoly2}
\frac{\det (\lambda - (A+\tau uv^\top  ))}{q(\lambda)}= m(\lambda) -\tau p_{uv}(\lambda) ,
\end{equation}
which finishes the proof of (i).

(ii) Assume that $\lambda_0\in\sigma(A+\tau_1 uv^\top)\cap\sigma(A+\tau_2 uv^\top)$ with
$\tau_1\neq\tau_2$. By (i) $\lambda_0$ is either a root of $q$, or a common root of the
polynomials $p_{uv,\tau_1}$ and $p_{uv,\tau_2}$. In the former case $\lambda_0$ clearly
belongs to $\sigma(A)$, in the latter case $\lambda_0$ is a root of $(\tau_1-\tau_2)p_{uv}$
and consequently of $m$. Hence, $\lambda_0\in\sigma(A)$ as well.

(iii) By Lemma \ref{pdegl-1}, for generic $u$ and $v$ and all
$\tau\in\mathbb{C}\setminus\set0$ the polynomials $p_{uv,\tau}$ and $m$ do not have common
roots and consequently $q$ is the greatest common divisor of the characteristic polynomials of
$A$ and $A+\tau uv^\top $. Hence, for  generic $u$ and $v$ the roots of $p_{uv,\tau}$ are
precisely the eigenvalues of $B(\tau )$ which are not eigenvalues of $A$. Since the $\deg
p_{uv,\tau}=\deg m$, there are exactly $\deg m$, counting algebraic multiplicities,
eigenvalues of $A+\tau uv^\top $ which are not eigenvalues of $A$.

\end{proof}

Note that by Lemma \ref{prop:1} for each $\tau\neq 0$ the eigenvalues in
$\sigma(B(\tau)\setminus\sigma(A)$ are non--derogatory. However, the proposition above does
not say, that for each $\tau\neq 0$ the eigenvalues in $\sigma(B(\tau)\setminus\sigma(A)$ are
simple. Obviously, for a fixed value of $\tau$ and generic $u$ and $v$ the eigenvalues in
$\sigma(B(\tau)\setminus\sigma(A)$ are simple, as follows from the H\"ormander--Mellin result,
but this is a weaker statement.

\section{The Jordan structure of $A+\tau uv^\top $ at the eigenvalues of $A$.}

The theorem below shows that the Jordan structure of $B(\tau)$ at the eigenvalues of $A$ is
constant for all $\tau\neq0$. The technique of the proof was used in  \cite{MMRR1} to reprove
the H\"ormander--Mellin result.

\begin{thm}\label{thm:generict}
Let $A\in\mathbb{C}^{n\times n}$ and let \eqref{JordanA}, \eqref{n's}  be the Jordan form of
$A$. Then for generic $u$ and $v$ and all $\tau \in\mathbb{C}\setminus\set0$ the sizes of the
Jordan blocks of $A+\tau uv^\top $ corresponding to the eigenvalue $\lambda_j$ are $n_{j,2}
\geq \dots \geq n_{j,k_j}$, for $j=1, \ldots , r$.
\end{thm}

\begin{proof}
Using the transformation \eqref{matrixtransform} we can assume that $A$ is in the Brunovsky
canonical form.  Denote by $e_{j,l}$ ($j=1,\dots ,r$, $l=1, \dots , n_{j,1}+n_{j,2}+ \cdots
+n_{j,k_j}$) the vector with one on the $l$--th position in the $j$--th block and zeros
elsewhere. Then the following sequences are Jordan chains of $A+\tau uv^\top $ corresponding
to the eigenvalue $\lambda_j$ ($j=1,\dots r$):
\begin{equation}\label{Jordanchains}
\begin{array}{l}
e_{j,1}-e_{j,n_{j,1}+1},\dots,e_{j,n_{j,2}}-e_{j,n_{j,1}+n_{j,2}};\\
e_{j,1}-e_{j,n_{j,1}+n_{j,2}+1},\dots,e_{j,n_{j,3}}-e_{j,n_{j,1}+n_{j,2}+n_{j,3}};\\
\vdots\\
e_{j,1}-e_{j,n_{j,1}+\cdots+n_{j,k_j-1}+1},\dots,e_{j,n_{j,k_j}}-e_{j,n_{j,1}+\cdots+n_{k_j-1}+n_{j,k_j}}.\\
\end{array}
\end{equation}
Hence, we see that for generic $u$ and $v$ there are Jordan chains of $A+\tau uv^\top $ of
lengths $n_{j,2} \geq \dots \geq n_{j,k_j}$ corresponding to the eigenvalue $\lambda_j$.
(Obviously, if $k_j=1$ then $\lambda_j$ is not an eigenvalue of $A+\tau uv^\top $). By
Proposition \ref{leigs} the dimension of the algebraic eigenspace corresponding to
$\sigma(B(\tau))\setminus\sigma(A)$ is $\deg m=n_{j,1}+\cdots +n_{r,1}$. Hence, none of the
Jordan chains in \eqref{Jordanchains} can be extended and the proof is finished.
\end{proof}


\section{The large $\tau$ asymptotic of eigenvalues of $B(\tau)$.}\label{s:larget}

In this  section it is shown that the eigenvalues of $B(\tau)$ that are not eigenvalues of $A$
tend with $\tau\to\infty$ to the roots of the polynomial $p_{uv}$, except one eigenvalue that
goes to infinity. This behavior is again generic in $u$ and $v$.

\begin{thm}\label{thm:largeasymp}
Let $A\in\mathbb{C}^{n\times n}$. Then for generic $u,v\in\mathbb{C}^n$ there exist
differentiable functions
 $$
 \mu_1,\dots,\mu_l:\set{\tau\in\mathbb{C}:|\tau|>\tau_0}\to \mathbb{C},
 $$
with $l=\deg m$ and some $\tau_0>0$,  such that
\begin{itemize}
\item[(i)] $\sigma(B(\tau))\setminus\sigma(A)=\set{\mu_j(\tau):j=1,\dots,l}$ for
$|\tau|>\tau_0$;
\item[(ii)] $\mu_j\neq\mu_{j'}$ for  $j,j'=1,\dots,l$, $j\neq j'$;
\item[(iii)] $\mu_1(\tau),\dots, \mu_{l-1}(\tau)$  tend  with $|\tau|\to\infty$ to the $l-1$ roots of the polynomial
$p_{uv}$;
\item[(iv)] $\mu_l(\tau)/\tau\to v^\top u$ with $\tau\to\infty$.
\end{itemize}
\end{thm}

The theorem says, in other words, that as $\tau$ goes to $\infty$ the eigenvalues of $B(\tau)$
which are not eigenvalues of $A$ are simple, exactly $l-1$ of them are approximate the roots
of $p_{uv}$ and one goes to infinity, asymptotically along the ray in the complex plane going
from zero through the number $v^\top  u$.

\begin{proof}
By Lemma \ref{pdegl-1} there are $l-1$ simple roots of the polynomil $p_{uv}$, let us denote
them by $\lambda_1,\dots,\lambda_{l-1}$.  Let $\varepsilon>0$ be such that the closed discs
 $$
 C_j(\eps)=\{\lambda\in\mathbb{C} : |\lambda -\lambda_j| \leq \varepsilon\},\qquad j=1,\dots l-1
 $$
do not intersect. Consider the polynomials
$$
{q}_\tau(\lambda)=\frac{1}{\tau}m(\lambda) -p_{uv}(\lambda),\quad \tau>\tau_0
$$
 and  observe that $\frac{1}{\tau}m(\lambda)$ converges with $|\tau|\to\infty$ uniformly to zero on $\bigcup_{j=1}^{l-1}C_j(\eps)$. By the Rouche theorem there is a $\tau_0>0$ so that for
$|\tau|>\tau_0$ the polynomial $q_\tau$ has exactly one simple root $\mu_j(\tau)$ in each of the
sets $C_j(\eps)$, $j=1,\dots,l-1$. Hence, the root $\mu_l(\tau)\notin
\bigcup_{j=1}^{l-1}C_j(\eps) $ is simple as well. By simplicity of the roots we get
$q_\tau'(\mu_j(\tau))\neq 0$ for $j=1,\dots l$, $|\tau|>\tau_0$. Hence, by implicit function
theorem the functions $\mu_1(\tau),\dots,\mu_l(\tau)$ are differentiable. Recalling that
$\sigma(B(\tau))\setminus\sigma(A)$ consists by Proposition \ref{leigs} precisely of the roots
of  $q_\tau(\lambda)$ finishes the proof of  (i) and (ii). Letting $\eps\to 0$ we obtain
(iii). To prove (iv) note that
 $$
\sigma\left(\frac1\tau B(\tau)\right)=\sigma\left(\frac1\tau A\right)\cup\set{\frac{\mu_1(\tau)}\tau,\dots,\frac{\mu_l(\tau)}\tau},\quad |\tau|>\tau_0.
$$
As $\tau\to\infty$ the matrix $\tau^{-1}B(\tau)$  converges  to the rank one matrix
 $uv^\top $ and thus $\mu_l(\tau)/\tau$ converges to $v^\top u$.
 \end{proof}

\begin{rem}\label{desc} In Figure
\ref{nodoublevlees} the roots of the polynomial $p_{uv}$ are marked with black circles, and
the asymptotic ray $y=(v^\top u)x$ is the dashed line.
\end{rem}

\section{Triple eigenvalues of $B(\tau)$.}
In this section we show that for generic $u,v$ there are no triple eigenvalues in
$\sigma(B(\tau))\setminus\sigma(A)$ for all $\tau\in\mathbb{C}$. In particular there are
generically no triple crossings of the eigenvalue curves.

\begin{thm}\label{thm:2}
Let $A\in\mathbb{C}^{n\times n}$. Then for generic $u,v\in\mathbb{C}^n$ and for all
$\tau\in\mathbb{C}$ the algebraic multiplicity of the eigenvalues of $A+\tau uv^\top$ that are
not eigenvalues of $ A$  is at most two.
\end{thm}

\begin{proof}
Suppose that $u$ and $v$ are such that for some $\tau\in\mathbb{C}$ the matrix $B(\tau)$ has
an eigenvalue $\lambda_0\notin\sigma(A)$ of multiplicity at least three. Then by Lemma
\ref{prop:1} $B(\tau)$ has a Jordan block of size at least three at $\lambda_0$. Consequently,
by Proposition \ref{leigs},  $\lambda_0$ is a triple root of $p_{uv,\tau}$, i.e.
\begin{align*}
 m(\lambda_0) -\tau p_{uv}(\lambda_0) &= 0 , \\
 m^\prime(\lambda_0) -\tau p_{uv}^\prime(\lambda_0)&=0, \\
 m^{\prime\prime}(\lambda_0) -\tau p_{uv}^{\prime\prime}(\lambda_0)&=0.
\end{align*}
Solving for $\tau$ from the first equation and substituting  in the second and third we obtain
\begin{align*}
& m^\prime(\lambda_0)p_{uv}(\lambda_0) -m(\lambda_0) p_{uv}^\prime(\lambda_0)=0, \\
& m^{\prime\prime}(\lambda_0)p_{uv}(\lambda_0) -m(\lambda_0)
p_{uv}^{\prime\prime}(\lambda_0)=0.
\end{align*}
Let $s$ be the greatest common divisor of $m$ and $m'$. 
Since $\lambda_0$ does not belong to
$\sigma(A)$, it is a common root of the polynomials
\begin{align}
f_{uv}=& \frac{m^\prime}{s}p_{uv} -\frac ms p_{uv}^\prime, \label{f}\\
g_{uv}=& m^{\prime\prime}p_{uv} -m\label{f'} p_{uv}^{\prime\prime}.
\end{align}
Therefore, $\det S(f_{uv},g_{uv})=0$. Summarizing, we showed so far that the set of all $u$
and $v$ for which there exists $\tau\in\mathbb{C}$ such that the matrix $B(\tau)$ has an
eigenvalue $\lambda_0\notin\sigma(A)$ of multiplicity at least three is contained in the set
of all $u,v\in\mathbb{C}^n$ such that $\det S(f_{uv},g_{uv})=0$. Clearly $\det
S(f_{uv},g_{uv})$ is a polynomial in the coordinates of $u$ and $v$. We show now, that it is a
nonzero polynomial, i.e. that for some $u,v$ the polynomials $f_{uv}$, $g_{uv}$ do not have a
common root, which will finish the proof. First consider the case  $\deg m=1$. Then for
generic $u,v$ the polynomial $p_{uv}$ is a constant nonzero polynomial and thus $f_{uv}$ is a
constant nonzero polynomial as well. Therefore, it does not have common roots with $g_{uv}$.
Now let us turn to the case $\deg m>1$.
 Observe that for every
$b\in\mathbb{C}$ there exist $u_b,v_b$ such that $p_{u_bv_b}(\lambda)=\lambda-b$. Then
\begin{align*}
&f_{u_bv_b}(\lambda)= \frac{m^\prime}{s}(\lambda)(\lambda-b) -\frac ms(\lambda), \\
&g_{u_bv_b}(\lambda)= m^{\prime\prime}(\lambda)(\lambda-b).
\end{align*}
Let $\mu_1,\dots \mu_{l-2}$ be the roots of $m''$. Note that $ \frac{m'}s(\mu_j)=0$ implies
$\frac ms(\mu_j)\neq0$ due to the definition of $s$. Therefore, one can  find
$b_0\in\mathbb{C}\setminus\sigma(A)$ such that
$$
\frac{m'}s(\mu_j)\cdot b_0\neq -\frac ms(\mu_j)-\frac{m'}s(\mu_j)\cdot\mu_j,\quad j=1,\dots
l-2.
$$
Consequently, $f_{u_{b_0}v_{b_0}}$ and $g_{u_{b_0}v_{b_0}}$ do not have a common root.
\end{proof}

Obviously, the  result  holds only generically. One can easily construct  a matrix $A$ and
vectors $u$ and $v$ such that $B(\tau)$ will have an  eigenvalue of a given multiplicity for a
given $\tau_0$. Namely, let $A_0=J_k(0)$  and let $u,v$ be any two  vectors for which
$A=A_0-\tau_0 uv^\top $ has $k$ different eigenvalues. Then $A+\tau uv^\top =J_k(0)$.

\section{Double eigenvalues of $B(\tau)$}\label{sec:doubleev}

\begin{thm}\label{prop:3}
Let  $A\in\mathbb{C}^{n\times n}$. Then  generic $u,v\in\mathbb{C}^n$ there are at most $2\deg
m-2$ values of the parameter $\tau\in\mathbb C$ for which there exists an eigenvalue of
$A+\tau uv^\top$ of multiplicity at least two, which is not an eigenvalue of $A$.
\end{thm}

\begin{proof}
Note that for all $\tau\in\mathbb R\setminus\{0\}$ the matrix $B(\tau)$ has a double
eigenvalue if and only if the polynomials $p_{uv,\tau}$ and $p'_{uv,\tau}$  have a common
zero, see Proposition \ref{leigs}. Write $m$ and $p_{uv}$ as
$$
m(\lambda)  = \lambda^l + \sum_{j=0}^{l-1} a_j \lambda^j, \qquad p_{uv}(\lambda) =
\sum_{j=0}^{l-1} p_{j} \lambda^j.
$$
Then the polynomials $p_{uv,\tau}$ and $p'_{uv,\tau}(\lambda)$ are given by
\begin{align*}
p_{uv,\tau}(\lambda)&= \lambda^l + \sum_{j=0}^{l-1} (a_j-\tau p_{j})\lambda^j ,
\\
p'_{uv,\tau}(\lambda)&=l\lambda^{l-1} +\sum_{j=1}^{l-1} j(a_j-\tau p_{j})\lambda^{j-1} .
\end{align*}
Consider the Sylvester resultant matrix $S(p_{uv,\tau},p'_{uv,\tau})\in\mathbb C^{(2l-1)\times
(2l-1)}$
 and let
 $$
 G(u,v,\tau)=\det S(p_{uv,\tau},p'_{uv,\tau}).
 $$
Then  $G(u,v,\tau)=0$  if and only if there is an eigenvalue of $B(\tau)$ of multiplicity at
least two, which is not an eigenvalue of $A$. Computing the determinant  $G(u,v,\tau)$ by
development of \eqref{Sylwester}  according to the first column (note that $a_{n_1}=1$,
$b_{n_2}=l$), one sees that it is the sum of constant in $\tau$ multiples of two determinants
of size $(2l-2)\times (2l-2)$, the entries of which are linear polynomials in $\tau$, or
constants. Using the fact that the determinant of a $k\times k$ matrix is a polynomial of
degree $k$ in the entries of the matrix, we see that $G(u,v,\cdot)$ is a polynomial of degree
at most $2l-2$ in the  variable $\tau$. This means that for any $A$, $u$ and $v$ the
polynomial $G(u,v,\cdot)$ has at most $2l-2$ zeros or is identically zero. However, by Theorem
\ref{thm:largeasymp}  we already know that for generic $u,v$ there exists $\tau_0\geq0$ such
that for $|\tau|>\tau_0$ the spectrum $\sigma(B(\tau))\setminus \sigma(A)$ consists of simple
eigenvalues only and consequently $G(u,v,\tau)\neq0$. Thus for generic $u,v$ the polynomial
$G(u,v,\cdot)$ has at most $2l-2$ roots and the theorem is proved.
\end{proof}

The last result of  this paper considers the real parameter $\tau$. Together with Proposition
\ref{leigs}(ii) it shows why the crossing of the eigenvalue curves in \eqref{problemstatement}
do not appear in numerical simulations, except possibly the crossings at $\sigma(A)$.

\begin{thm}\label{propC}
Let $A\in\mathbb{C}^{n\times n}$ and let $V$  be the set of all pairs
$(u,v)\in\mathbb{C}^{2n}$ for which there exists $\tau\in\mathbb{R}$ such that $A+\tau uv^\top
$ has a double eigenvalue, which is not an eigenvalue of $A$. Then $V$ is closed, with empty
interior and has the $4n$--dimensional Lebesgue measure zero.
\end{thm}

\begin{proof}
  As in the proof of Theorem \ref{prop:3} we note that
 $$
 V=\set{(u,v)\in\mathbb{C}^{2n}:\exists_{\tau\in\mathbb{R}\setminus\set0}\ G(u,v,\tau)=0}.
 $$
Since the zeros of a polynomial depend continuously on its coefficients, the set is
$\mathbb{C}\setminus V$ is open. To prove that $V$  is of $4n$--dimensional Lebesque measure
zero (and consequently has an empty interior) consider  the set
$$
U_0:=\set{(u,v)\in\mathbb{C}^{2n}: \exists_{ \lambda\in\mathbb{C}}\ f_{uv}(\lambda)=f'_{uv}(\lambda)=0},
$$
where $f_{u v}$ is defined as in \eqref{f}. Note that
$$
U_0=\set{(u,v)\in\mathbb{C}^{2n}: \exists_{ \lambda\in\mathbb{C}}\
f_{uv}(\lambda)=g_{uv}(\lambda)=0},
$$
where $g_{uv}$ is defined as in \eqref{f'}. Indeed, this follows from
$$
s^2 f'_{uv}=sg_{uv}-s'f_{uv}
$$
and from the fact that the polynomials $s$ and $f_{uv}$ do not have common roots. Hence, it
follows from the proof of Theorem \ref{prop:3} that the set $U_0$ is a proper algebraic subset
of $\mathbb{C}^{2n}$.

 Recall that by Lemma \ref{pdegl-1} the set
$$
 U_1=\set{(u,v)\in\mathbb{C}^{2n}:\deg p_{uv}<l-1},
$$
is also a proper algebraic subset of $\mathbb{C}^n$. Observe that for $(u,v)\notin U_1$ one
has $\deg f_{uv}=k$, where $k:=\max\set{(r-1)(l-1),r(l-2)}$,  $l=\deg m$ and $r=\deg\frac ms$
is the number of eigenvalues of $A$. To see this let $(u,v)\notin U_1$. In the case
$(r-1)(l-1)\neq r(l-2)$ it is clear that $\deg f_{uv}=k$. In the case when $(r-1)(l-1)=r(l-2)$
note that although the degrees of both summands in \eqref{f} coincide, the leading coefficient
does not cancel. Indeed, the leading coefficients of $\frac{m'}s p_{uv}$ and $\frac ms
p'_{uv}$ are respectively $l\alpha$ and $(l-1)\alpha$, where $\alpha$ is the leading
coefficient of $p_{uv}$.

Consequently,
$$
V_0:=\mathbb{C}^{2n}\setminus (U_0\cup U_1)
$$
 is an open and nonempty set. Note that for each
$(u,v)\in V_0$ the function $f_{uv}$ has precisely $k$ zeros $\lambda_1(u,v),\dots,
\lambda_k(u,v)$ and they are all not in $\sigma(A)$.  Since $f_{uv}(\lambda_j(u,v))=0$ and
$(u,v)\notin U_0$, one has $f'_{uv}(\lambda_j(u,v))\neq 0$, $j=1,\dots,k$. Therefore, by the
implicit function theorem, the functions $\lambda_1(u,v),\dots, \lambda_k(u,v)$ can be chosen
as holomorphic functions on $V_0$.  Note that
$$
V\subseteq U_0\cup U_1\cup\bigcup_{j=1}^k V_j,
$$
with
$$
V_j= \set{u\in V_0: \exists_{\tau\in\mathbb{R}\setminus\set0}\ m(\lambda_j(u,v))- \tau
p_{uv}(\lambda_j(u,v))=0}
$$ $$
=\set{u\in V_0: \frac{p_{uv}(\lambda_j(u,v))}{m(\lambda_j(u,v))}\in\mathbb{R}},\quad
 j=1,\dots, k.
$$
Observe that  the functions
 $$
 V_0\ni (u,v)\mapsto \frac{p_{uv}(\lambda_j(u,v))}{m(\lambda_j(u,v))}=v^\top(\lambda_j(u,v)-A)^{-1}u \in\mathbb{C},\quad
j=1,\dots,j
 $$
are holomorphic and nonconstant on every connected component of $V_0$. By the uniqueness
principle each of the  sets $V_j$ ($ j=1,\dots, k$) is of $4n$--dimensional Lebesgue measure
zero. Hence, $V_0$, and in consequence $V$ as well, are of $4n$--dimensional Lebesgue measure
zero.

\end{proof}

In the infinite dimensional case the function $Q(z)=-\seq{(\lambda-A)^{-1}u,u}$ is a very useful tool for studying spectra of one dimensional perturbations of selfadjoint operators, or even more generally, spectra of finite dimensional selfadjoint extensions of symmetric operators. The key point is solving the equation $Q(z)=-1/\tau$ and as it can be seen this technique was a motivation for the proof above. This approach can be found e.g. in  \cite{HSSW} in the Hilbert space context and in \cite{DHS1,DHS3,SWW} in the Pontryagin space setting.

\end{document}